\documentclass[12pt]{amsart}
\usepackage{amssymb}
\usepackage{eucal}
\newcommand{\rr}{\mathbb{R}}
\newcommand{\p}{\partial}
\newcommand{\zz}{\mathbb{Z}}

\newtheorem{theorem}{Theorem}[section]

\author{A. Alexandrou Himonas, Gerard Misio\l ek, Gustavo Ponce, and Yong Zhou}
\title[Propagation speed for the C-H equation] 
{Persistence Properties and Unique Continuation of 
solutions of the Camassa-Holm equation} 

\begin{document}

\maketitle

\begin{abstract}  It is shown that a strong solution  of the Camassa-Holm equation, 
initially decaying exponentially  together with its spacial derivative,
must be  identically equal to zero if it also decays exponentially
at a later time. In particular, a strong solution of the Cauchy problem
with compact initial profile  can not  be compactly supported at any
later time unless it is the zero solution. 
\end{abstract}

\numberwithin{equation}{section} 

\section{Introduction}

This work is concerned with the nonperiodic Camassa-Holm equation 
\begin{equation}
\label{1.1} 
\p_t u-\p_{t}\p^2_{x}u+3u\p_x u - 2\p_x u \p^2_x u - u\p_x^3 u = 0,
\;\;\;\;t,\,x\in\rr. 
\end{equation} 
This equation appears in the context of hereditary
symmetries studied by Fuchssteiner and Fokas \cite{ff}.
It was first written explicitly, and derived physically as a water wave equation
by Camassa and Holm \cite{ch}, who also studied its solutions.
Equation (1.1) is remarkable for its properties such as infinitely
many conserved integrals, bi-hamiltonian structure or its non-smooth
travelling wave solutions known as ``peakons" (see formula (1.9)).
It was also derived as an equation for geodesics of the $H^1$-metric
on the diffeomorphism group, see \cite{mi}.
For a discussion of how it relates to the theory of hereditary
symmetries see \cite{f}. 
The inverse scattering approach to the Camassa-Holm equation has also
been developed in several works,
for example see \cite{ch}, \cite{comc}, \cite{mc}, \cite{bss}, and references therein.

A considerable amount of work has been devoted to the study 
of the corresponding Cauchy problem in both nonperiodic and 
periodic cases. 
Among these results, of relevance to the present paper will be 
the fact that (\ref{1.1}) is locally well-posed (in Hadamard's sense) 
in $H^s(\mathbb{R})$ for any $s>3/2$, see for example \cite{lo}, 
\cite{rb}, \cite{d}. The long time behaviour of  solutions has been 
studied and conditions which guarantee their global existence and
their finite blow up have been deduced. In particular,  in \cite{mc}  a necessary and sufficient
condition was established on the initial datum to guarantee finite time singularity
formation for the corresponding strong solution. For further results in this direction we refer to
\cite{mc}, \cite{CE}  and 
 the survey article  \cite{mo} and references therein.
For well-posedness results in the periodic case we refer to 
\cite{hm1}, \cite{mi}, and \cite{dkt}, where the equation is studied in its 
integral-differential form (see \eqref{1.2} below) as an ODE 
on the space of diffeomorphisms of the circle. A recent 
result  demonstrating that the solution map $u_0\rightarrow u$
for the Camassa-Holm equation is not locally uniformly continuous 
in   Sobolev spaces  can be found in \cite{hm2}. 

Also the Camassa-Holm equation has been studied as an integrable infinite-dimensional
Hamiltonian system, and several works have been devoted to several aspect of its
scattering setting, see \cite{ch}, \cite{comc}, \cite{mc}, \cite{bss} and references therein.

It is convenient to rewrite the equation in 
its formally equivalent integral-differential form
\begin{equation}
\label{1.2}
\p_t u + u\p_x u +\p_x G\ast \left( u^2+\frac{1}{2}(\p_x u)^2 \right) =0,
\end{equation}
where $G(x)=e^{-|x|}$.

Our first objective here is to formulate decay conditions on a solution, 
at two distinct times, which guarantee that $u \equiv 0$ 
is the unique solution of equation (\ref{1.1}). 
The idea of proving unique continuation results for nonlinear dispersive 
equations under decay assumptions of the solution at two different times 
was motivated by the recent works \cite{EKPV1}, \cite{EKPV2} on 
the nonlinear Schr\"odinger and the $k$-generalized Korteweg-de Vries 
equations respectively. 

In  the recent works \cite{C}, \cite{He} and \cite{Z} it was shown 
that $u$ cannot preserve compact support in a non-trivial time interval 
( i.e. for $\,t \in [0,\epsilon],\;\epsilon>0$) 
except if $u\equiv 0$. However, this result does not preclude the possibility
of the solution  having compact support at a later time.
In fact, in \cite{Z} the question concerning the possibility of 
a smooth solution of \eqref{1.1}  having compact support 
at two different times was explicitly stated. In particular, our first result,
Theorem \ref{Theorem 1}, gives a negative answer to this question.

\begin{theorem}
\label{Theorem 1} 
Assume that for some $T>0$ and $s>3/2$
\begin{equation}
\label{1.3}
u\in C([0,T]\,:\,H^s(\rr)) 
\end{equation}
is a strong solution of the IVP associated to the equation \eqref{1.2}. 
If 
$u_0(x)=u(x,0)$
satisfies that for some $\alpha\in (1/2, 1)$
\begin{equation}
\label{1.4}
|u_0(x)|\sim o(e^{-x}),\;\;\;\;\text{and}\;\;\;\;|\p_x u_0(x)|\sim
O(e^{-\alpha x})\;\;\;\;\text{as}\;\;x\uparrow \infty,
\end{equation}
and there exists $t_1\in (0,T]$ such that 
\begin{equation}
\label{1.5}
|u(x,t_1)|\sim o(e^{-x})\;\;\;\;\text{as}\;\;\;x\uparrow \infty,
\end{equation}
then $u\equiv 0$.
\end{theorem}

\underline {Remarks} (a)  Theorem \ref{Theorem 1}  holds with the corresponding  
decay hypothesis in \eqref{1.4}-\eqref{1.5}  stated for $x<0$.

(b) The time interval $[0,T]$ is the maximal existence time 
interval of the strong solution.
This guarantees that the solution is uniformly bounded in 
the $H^s$-norm in this interval (see \eqref{1.22}), and that 
our solution is the strong limit of smooth ones such that 
the integration by parts in the proof (see \eqref{1.31}, \eqref{1.39}) 
can be justified.

\vskip 0.3cm 

The proof of Theorem \ref{Theorem 1} will be a consequence of 
the following result concerning some persistence properties 
of the solution of the equation \eqref{1.2}
in $L^{\infty}$-spaces with exponential weights.

\begin{theorem}\label{Theorem 2} 
Assume that for some $T>0$ and $s>3/2$
\begin{equation}
\label{1.6}
u\in C([0,T]\,:\,H^s(\rr)) 
\end{equation}
is a strong solution of the IVP associated to 
the equation \eqref{1.2} and that
$u_0(x)=u(x,0)$ satisfies that for some $\theta \in (0,1)$
\begin{equation}
\label{1.7}
|u_0(x)|,\;\;\;|\p_x u_0(x)|\,\sim O(e^{-\theta
x})\;\;\;\text{as}\;\;\;\;x\uparrow \infty.
\end{equation}
Then
\begin{equation}
\label{1.8}
|u(x,t)|,\;\;\;|\p_x u(x,t)|\,\sim O(e^{-\theta
x})\;\;\;\text{as}\;\;\;\;x\uparrow \infty,
\end{equation}
uniformly in the time interval $[0,T]$.
\end{theorem}

The following result establishes the optimality of 
Theorem \ref{Theorem 1} and tells us that a strong non-trivial solution 
of \eqref{1.2} corresponding to data with fast decay at infinity 
will immediately behave asymptotically, in the  $x$-variable at infinity,  
as the \lq\lq peakon" solution
\begin{equation} \label{1.20.b}
v_c(x,t)=c\,e^{-|x-ct|},\;\;\;\;t>0.
\end{equation}

\begin{theorem}\label{Theorem 3} 
Assume that for some $T>0$ and $s>3/2$
\begin{equation}
\label{1.9}
u\in C([0,T]\,:\,H^s(\rr)) 
\end{equation}
is a strong solution of the IVP associated to the equation \eqref{1.2} 
and that
$u_0(x)=u(x,0)$ satisfies that for some $\alpha\in(1/2,1)$
\begin{equation}
\label{1.10}
|u_0(x)| \sim  O(e^{-x}),\;\;\;
|\p_x u_0(x)|\,\sim O(e^{-\alpha x})\;\;\;
\text{as} \;\; x \uparrow \infty 
\end{equation}
for some $\alpha \in (1/2,1)$. 
Then
\begin{equation}
\label{1.11}
|u(x,t)|\sim O(e^{- x})\;\;\;\text{as}\;\;x\uparrow \infty,
\end{equation}
uniformly in the time interval $[0,T]$.
\end{theorem}

 In the case when the solution $u(x,t)$ possesses  further regularity
and its data $u_0$ has stronger decay properties we shall  give a more precise description of its  
behavior at infinity in the space variable.
 
As it was noted in both \cite{C} and \cite{Z} in the case of compactly 
supported initial data $u_0$ the difference $h(x,t)$ of the solution and its second derivative, i.e. 
\begin{equation}
\label{1.11a}
h(x,t) = (1 -\p_x^2 )u(x,t),
\end{equation}
remains compactly supported. Thus, if $u_0$ is supported in the interval $[a,b]$ in its life-span, one has that $ h(x,t)$ has compact support in the time interval
$[\eta(a,t), \,\eta(b,t)]$, (for the definition of $\eta(\cdot,\cdot)$ see \eqref{1.48}).

\begin{theorem} \label{Theorem 4} 
Assume that for some $T>0$ and $s>5/2$ 
\begin{equation} \label{1.12c}
u\in C([0,T]\,:\,H^s(\rr)) 
\end{equation}
is a strong solution of the IVP associated to the equation \eqref{1.2}.
 
(a) If $u_0(x) =u(x,0)$ has compact support, 
then for any $t \in (0,T]$
\begin{equation} \label{1.13c} 
u(x,t)=\begin{aligned}
\begin{cases}
& c_+(t)\, e^{-x},\;  \;\;\;\text{for}\;\;\;\; x> \eta(b,t),\\
&c_{-}(t) \,e^{x},\;  \;\;\;\;\;\text{for}\;\;\;\; x< \eta(a,t).
\end{cases}
\end{aligned}
\end{equation} 

(b) If for some $\mu >0$ 
\begin{equation}
\label{1.13d}
\p_x^j u_0\sim  O(e^{-(1+\mu) |x|})\;\;\;\;\;\text{as}\;\;\;\;\;|x|\uparrow \infty\;\;\;\;\;j=0,1,2,
\end{equation}
then for any $t \in (0,T]$
\begin{equation}
\label{1.13e}
h(x,t) = (1 -\p_x^2 )\,u(x,t)\sim  O(e^{-(1+\mu) |x|})\;\;\;\;\;\;\text{as}\;\;\;\;\;|x|\uparrow \infty,
\end{equation}
and
\begin{equation}
\label{1.13f}
\lim_{x\to \pm \infty} e^{\pm x} u(x,t)=c_{\pm}(t),
\end{equation}
where in \eqref{1.13c}, \eqref{1.13f}  $ c_{+}(\cdot),\;c_{-}(\cdot)$ denote continuous non-vanishing functions,
with $\,c_{+}(t)>0$ and $\,c_{-}(t)<0$ for $t\in(0, T]$.
Furthemore, $ c_{+}(\cdot) $ is strictly increasing function, while $ c_{-}(\cdot) $ is strictly decreasing.

\end{theorem} 

Theorem \ref{Theorem 4}  tells us that, as long as it exists, the solution
$u(x,t)$ is positive at infinity and negative at minus infinity
regardless of the profile of a fast-decaying data $u_0 \neq 0$.

\vskip.2in


\section{Proof of the results}
First, assuming the result in Theorem \ref{Theorem 2} 
we shall prove Theorem \ref{Theorem 1}. 

\begin{proof}[Proof of Theorem \ref{Theorem 1}]
\label{Proof of Theorem 1}
Integrating equation \eqref{1.2} over the time interval 
$[0,t_1]$ we get
\begin{equation}
\label{1.12}
u(x,t_1)-u(x,0)+\int_0^{t_1}u\p_x u(x,\tau) d\tau +\int_0^{t_1}
\p_xG\ast (u^2+\frac{1}{2}(\p_xu)^2)(x,\tau)d\tau=0.
\end{equation}
By hypothesis \eqref{1.4} and \eqref{1.5} we have 
\begin{equation}
\label{1.13}
u(x,t_1)-u(x,0) \sim o(e^{-x})\;\;\;\text{as}\;\;x\uparrow \infty.
\end{equation}
From \eqref{1.4} and Theorem \ref{Theorem 2} it follows that
\begin{equation}
\label{1.14}
\int_0^{t_1}u\p_x u(x,\tau) d\tau 
\sim  
O(e^{-2 \alpha x})\;\;\;\text{as}\;\;x\uparrow \infty,
\end{equation}
and so 
\begin{equation}
\label{1.14.b}
\int_0^{t_1}u\p_x u(x,\tau) d\tau  
\sim
o(e^{-x})\;\;\;\text{as}\;\;x\uparrow \infty.
\end{equation}

We shall show that if $u\ne 0$, then the last term in \eqref{1.12} 
is $O(e^{-x})$ but not $o(e^{-x})$.
Thus, we have 
\begin{equation} \label{1.15} 
\begin{split} 
\int_0^{t_1}
\p_xG\ast (u^2+\frac{1}{2}(\p_xu)^2)(x,\tau) \, d\tau 
&=\p_x G\ast \int_0^{t_1}
(u^2+\frac{1}{2}(\p_xu)^2)(x,\tau)\, d\tau  \\
&= 
\p_x G\ast \rho (x),
\end{split} 
\end{equation}
where by \eqref{1.4} and Theorem \ref{Theorem 2} 
\begin{equation}
\label{1.16}
0\leq \rho(x)\sim O(e^{-2\alpha x}),
\;\;\;\;\text{so that}\;\;\;\;  
\rho(x)\sim o(e^{-x})\;\;\;
\text{as}\;\;x\uparrow \infty.
\end{equation}
Therefore
\begin{equation}
\begin{split}
\p_xG\ast \rho (x) 
&= 
- \int_{-\infty}^{\infty} sgn(x-y)\,e^{-|x-y|}\,\rho(y) \,dy\\
&= 
-e^{-x}\,\int_{-\infty}^x e^y \rho(y) \, dy 
+ e^x\, \int_x^{\infty}e^{-y}\rho(y) \, dy.
\label{1.17}
\end{split}
\end{equation}
From \eqref{1.16} it follows that
\begin{equation}
\label{1.18}
e^x\, \int_x^{\infty}e^{-y}\rho(y)\, dy 
= 
o(1) e^x\, \int_x^{\infty}e^{-2y} \, dy
\sim o(1) e^{-x}\sim o(e^{-x}),
\end{equation}
and if $\,\rho \ne 0$ one has that
\begin{equation}
\label{1.19}
\int_{-\infty}^x\,e^y\rho(y) \,dy 
\geq 
c_0, \qquad \text{for}\quad x \gg 1.
\end{equation}
Hence, the last term in \eqref{1.15} and \eqref{1.17} 
satisfies 
\begin{equation}
\label{1.20}
- \p_x G\ast \rho(x) 
\geq 
\frac{c_0}{2} \,\,e^{-x}, \qquad \text{for} \quad x\gg 1
\end{equation}
which combined with \eqref{1.12}-\eqref{1.14} yields a contradiction. 
Thus, $\rho(x)\equiv 0$ and consequently $u\equiv 0$, 
see \eqref{1.15}. 
\end{proof} 
  
\medskip

\begin{proof}[Proof of Theorem \ref{Theorem 3}] \label{Proof of Theorem 3} 
This proof is similar to that given for Theorem \ref{Theorem 1} 
and therefore it will be omitted.
\end{proof}


We proceed to prove Theorem \ref{Theorem 2}. 

\begin{proof}[Proof of Theorem \ref{Theorem 2}] 
\label{Proof of Theorem 2}
We introduce the following notations
\begin{equation}
\label{1.21}
F(u) = u^2 +\frac{\,(\p_x u)^2}{2},
\end{equation}
and
\begin{equation}
\label{1.22}
M= \sup_{t\in[0,T]}\;\|u(t)\|_{H^s}.
\end{equation}

Multiplying the equation \eqref{1.2} by $u^{2p-1}$ with $p\in\zz^+$ 
and integrating the result in the $x$-variable one gets 
\begin{equation} \label{1.23}
\int_{-\infty}^{\infty} u^{2p-1}\p_tu \, dx 
+ 
\int_{-\infty}^{\infty} u^{2p-1} \,u\p_xu \, dx 
+ 
\int_{-\infty}^{\infty}
u^{2p-1}\,\p_x G\ast F(u) \,dx =0.
\end{equation}
The estimates
\begin{equation} \label{1.24}
\int_{-\infty}^{\infty} u^{2p-1}\p_tu \, dx 
= 
\frac{1}{2p}\, \frac{d\;}{dt}\|u(t)\|_{2p}^{2p} 
= 
\|u(t)\|_{2p}^{2p-1} \, \frac{d\;}{dt}\|u(t)\|_{2p} 
\end{equation}
and 
\begin{equation}
\label{1.25}
\left| \int_{-\infty}^{\infty} u^{2p-1} \,u\p_xu \, dx \right| 
\leq 
\|\p_x u(t)\|_{\infty}\,\|u(t)\|_{2p}^{2p} 
\end{equation}
and H\"older's inequality in \eqref{1.23} yield
\begin{equation} \label{1.26}
\frac{d\;}{dt}\|u(t)\|_{2p} 
\leq 
\|\p_x u(t)\|_{\infty} \|u(t)\|_{2p} + \|\p_x G\ast F(u)(t)\|_{2p} 
\end{equation}
and therefore, by Gronwall's inequality 
\begin{equation} \label{1.27}
\|u(t)\|_{2p} 
\leq 
\left( \|u(0)\|_{2p} 
+ \| \p_xG\ast F(u)(\tau)\|_{2p}\,d\tau \right) e^{Mt}.
\end{equation}
Since $f\in L^2 (\rr) \cap L^{\infty}(\rr)$ implies 
\begin{equation}
\label{1.28}
\lim_{q \uparrow \infty} \|f\|_{q} = \|f\|_{\infty}, 
\end{equation}
taking the limits in \eqref{1.27} 
(notice that $\p_x G\in L^1$ and $F(u)\in L^1\cap L^{\infty}$) 
from \eqref{1.28} we get 
\begin{equation} \label{1.29}
\|u(t)\|_{\infty} 
\leq 
\left( \|u(0)\|_{\infty} 
+ 
\int_0^t\,\|\p_xG\ast F(u)(\tau)\|_{\infty}\,d\tau \right) e^{Mt}.
\end{equation}

Next, differentiating \eqref{1.2} in the $x$-variable produces 
the equation
\begin{equation} \label{1.30}
\p_t\p_x u + u\p^2_x u +(\p_xu)^2  
+ 
\p^2_x G\ast \left( u^2+\frac{1}{2}(\p_x u)^2 \right) =0.
\end{equation}
Again, multiplying the equation \eqref{1.30} by $\partial_x u^{2p-1}$ 
($p\in\zz^+$) integrating the result in the $x$-variable 
and using integration by parts 
\begin{equation}
\label{1.31}
\int_{-\infty}^{\infty} u\p^2_x u (\p_xu)^{2p-1} dx 
= 
\int_{-\infty}^{\infty} 
u \, \p_x\Big(\frac{(\p_xu)^{2p} }{2p}\Big) dx 
=
-\frac{1}{2p}\int_{-\infty}^{\infty} \p_x u (\p_xu)^{2p} \, dx 
\end{equation} 
one gets the inequality
\begin{equation}
\label{1.32}
\frac{d\;}{dt}\|\p_x u(t)\|_{2p} 
\leq 
2 \|\p_x u(t)\|_{\infty} \|\p_x u(t)\|_{2p} 
+
\|\p^2_xG\ast F(u)(t)\|_{2p} 
\end{equation} 
and therefore as before 
\begin{equation} \label{1.33}
\|\p_x u(t)\|_{2p} 
\leq 
\left( \|\p_x u(0)\|_{2p} 
+ 
\int_0^t\,\|\p^2_xG\ast F(u)(\tau)\|_{2p} \,d\tau \right) e^{2Mt}.
\end{equation}
Since $\p_x^2 G=G - \delta$, we can use \eqref{1.28} 
and pass to the limit in \eqref{1.33} to obtain
\begin{equation}
\label{1.34}
\|\p_x u(t)\|_{\infty} 
\leq 
\left( \|\p_xu(0)\|_{\infty} 
+ 
\int_0^t\,\|\p^2_xG\ast F(u)(\tau)\|_{\infty} \,d\tau \right) e^{2Mt}.
\end{equation}

We shall now repeat the above arguments using the weight 
\begin{equation}
\varphi_N(x)=\begin{cases}
1,\;\;\;\;\;\;\;\;x\leq 0,\\
e^{\theta x},\;\;\;\;\;x\in(0, N),\\
e^{\theta N},\;\;\;\;x\geq N 
\end{cases}
\label{1.35}
\end{equation}
where $N\in\zz^+$. Observe that for all $N$ we have 
\begin{equation} \label{1.36} 
0\leq \varphi'_N(x) 
\leq 
\varphi_N(x) \qquad \text{a.e.}\;\;\;x\in\rr.
\end{equation} 
Using notation in \eqref{1.21}, from equation \eqref{1.2} we obtain 
\begin{equation} \label{1.37}
\p_t (u\,\varphi_N) + (u\,\varphi_N) \p_x u +\varphi_N \,\p_x G\ast F(u)=0,
\end{equation}
while from \eqref{1.30} we get 
\begin{equation} \label{1.38}
\p_t (\p_x u \, \varphi_N) + u \, \p_x^2u \, \varphi_N 
+ 
(\p_xu \, \varphi_N) \p_x u + \varphi_N \,\p^2_x G\ast F(u)=0,
\end{equation}
We need to eliminate the second derivatives in the second term 
in \eqref{1.38}. 
Thus, combining integration by parts and \eqref{1.36} we find 
\begin{equation}
\begin{split} 
\Big| \int_{-\infty}^{\infty} u &\p_x^2 u \, \varphi_N(x) 
(\p_x u \,\varphi_N(x) )^{2p-1} dx\Big| \\
&= 
\Big| \int_{-\infty}^{\infty} u (\p_x u \, \varphi_N(x))^{2p-1} 
\left( \p_x(\p_xu \, \varphi_N(x) \right)  
- 
\p_x u \, \varphi'_N(x)) dx \Big| \\ 
&= 
\Big| \int_{-\infty}^{\infty} 
u \, \p_x \left( \frac{(\p_x u \, \varphi_N(x))^{2p}}{2p}\right) dx 
- 
\int_{-\infty}^{\infty} u \, \p_x u \, \varphi'_N(x) 
(\p_x u \,\varphi_N(x))^{2p-1} dx\Big|  \\
&\leq 
2\big( \|u(t)\|_{\infty}+\|\p_x u(t)\|_{\infty} \big) 
\|\p_x u\,\varphi_N\|_{2p}^{2p}
\end{split}
\label{1.39}
\end{equation}
Hence, as in the weightless case  \eqref{1.29} and  \eqref{1.34}, 
we get
\begin{equation}
\begin{split}
\|u(t)\varphi_N\|_{\infty} &+ \|\p_x u (t) \varphi_N\|_{\infty}
\leq  e^{2Mt} \left(\|u(0)\varphi_N\|_{\infty} 
+ 
\|\p_xu(0) \varphi_N\|_{\infty} \right)  \\
&+  
\,e^{2Mt}  
\int_0^t\, \left( \|\varphi_N \p_x G\ast F(u)(\tau)\|_{\infty} 
+ 
\|\varphi_N \p^2_xG\ast F(u)(\tau)\|_{\infty} \right)\,d \tau.
\end{split}
\label{1.40}
\end{equation}
A simple calculation shows that there exists $c_0>0$, depending only 
on $\theta \in (0,1)$  (see \eqref {1.7} and \eqref{1.35})  such that 
for any $N \in \zz^+$
\begin{equation}
\label{1.41}
\varphi_N(x) 
\int_{-\infty}^{\infty}\,e^{-|x-y|}\;\frac{1}{\varphi_N(y)} \,dy
\leq c_0.
\end{equation}

Thus, for any appropriate  function $f$ one sees that
\begin{equation}
\begin{split}
\big| \varphi_N\, \p_x G\ast f^2(x) \big| 
&= 
\left| \varphi_N(x) \int_{-\infty}^{\infty} 
sgn(x-y)\,e^{-|x-y|} f^2(y)\, dy \right|   \\
&\leq 
\varphi_N(x)\;\int_{-\infty}^{\infty} e^{-|x-y|} \,
\frac{1}{\varphi_N(y)} \varphi_N(y) f(y) f(y) \, dy   \\
&\leq 
\left( \varphi_N(x) \int_{-\infty}^{\infty} 
e^{-|x-y|}\;\frac{1}{\varphi_N(y)}\,dy \right) 
\|\varphi_N f\|_{\infty}\,\|f\|_{\infty} \\
&\leq 
c_0 \|\varphi_N f\|_{\infty}\,\|f\|_{\infty}.
\end{split}
\label{1.42}
\end{equation}
Since $\p_x^2 G=G - \delta$ the argument in \eqref{1.42} also shows that
\begin{equation}
\label{1.43}
\left| \varphi_N \,\p^2_x G\ast f^2(x) \right| 
\leq 
c_0 \|\varphi_N f\|_{\infty}\,\|f\|_{\infty}.
\end{equation}
Thus, inserting \eqref{1.42}-\eqref{1.43} into \eqref{1.40} 
and using \eqref{1.21}-\eqref{1.22} it follows that there exits 
a constant $\tilde c=\tilde c(M;T)>0$ such that
\begin{equation}
\begin{split}
\|u(t&) \varphi_N\|_{\infty} + \|\p_x u (t) \varphi_N\|_{\infty}
\leq 
\tilde c\big(\|u(0)\varphi_N\|_{\infty} +\|\p_xu(0) \varphi_N\|_{\infty}\big) \\
&+ 
\tilde c \int_0^t 
\big( \|u(\tau)\|_{\infty}+\|\p_x u(\tau)\|_{\infty} \big) 
\big( \|u(\tau)\varphi_N\|_{\infty} 
+ 
\|\p_x u(\tau)\varphi_N\|_{\infty} \big) d\tau \\
&\leq 
\tilde c \left( \|u(0)\varphi_N\|_{\infty} + \|\p_xu(0) \varphi_N\|_{\infty} 
+ 
\int_0^t \big( \|u(\tau)\varphi_N\|_{\infty} 
+ 
\|\p_x u(\tau)\varphi_N\|_{\infty}\big) d \tau \right).
\end{split}
\label{1.44}
\end{equation}
Hence, for any $N\in\zz^+$ and any $t\in [0,T]$ we have 
\begin{equation}
\begin{split}
\|u(t)\varphi_N\|_{\infty} + \|\p_x u (t) \varphi_N\|_{\infty}
&\leq 
\tilde c 
\big( \|u(0)\varphi_N\|_{\infty} + \|\p_xu(0) \varphi_N\|_{\infty}\big) \\
&\leq 
\tilde c 
\big(\|u(0) e^{\theta x}\|_{\infty} +\|\p_xu(0) e^{\theta x}\|_{\infty}\big).
\end{split}
\label{1.45}
\end{equation}
Finally, taking the limit as $N$ goes to infinity in \eqref{1.45} we find that 
for any $N\in\zz^+$ and any $t\in [0,T]$ 
\begin{equation}
\label{1.46}
\sup_{t \in[0,T]} 
\left( \|u(t)e^{\theta x} \|_{\infty} 
+ 
\|\p_x u (t) e^{\theta x}\|_{\infty} \right) 
\leq  
\tilde c 
\left( \|u(0) e^{\theta x} \|_{\infty} 
+ 
\|\p_xu(0) e^{\theta x}\|_{\infty} \right),
\end{equation}
which completes the proof of Theorem 2.
\end{proof}


It remains to prove  Theorem  \ref{Theorem 4}. 

\begin{proof}[Proof of Theorem \ref{Theorem 4}] \label{Proof of Theorem 4} 

A simple calculation shows that the solution $u$ 
of equation \eqref{1.1} satisfies the identity 
\begin{equation} \label{1.47} 
(1 - \p_x^2)u \circ \eta \, (\p_x \eta)^2  
= 
(1-\p_x^2)u_0, 
\end{equation} 
(it has a mechanical interpretation as conservation of spacial 
angular momentum). 
Here $\eta = \eta(x,t)$ is the flow of $u$, that is 
\begin{equation} \label{1.48} 
\begin{split} 
\begin{cases}
\dfrac{d\eta(x,t)}{dt} &= u(\eta(x,t),t),  \\ 
\eta(x,0) &= x, 
\end{cases}
\end{split} 
\end{equation} 
so that by the assumption and the standard ODE theory 
$t \to \eta(t)$ is a smooth curve of $C^1$-diffeomorphisms 
of the line, close to the identity map and defined on 
the same time interval as $u$ (see \cite{mi} for details 
in the periodic case). 
From \eqref{1.47} we then have 
\begin{equation} \label{1.49} 
u(x,t)  =  \int_{-\infty}^{\infty} e^{-|x-y|} h(y,t) dy 
= 
e^{-x} \int_{-\infty}^x e^y h(y,t) dy 
+ 
e^x \int_x^{\infty} e^{-y} h(y,t) dy, 
\end{equation} 
where 
\begin{equation} \label{1.50} 
h(x,t) =(1-\p_x^2)u(x,t) 
= 
\frac{(1-\p_x^2)u_0(\eta^{-1}(x,t))}{\big(\p_x\eta(\eta^{-1}(x,t),t)\big)^2}. 
\end{equation} 

Let us first prove part (a). 

Thus, from \eqref{1.50} it follows that if $u_0$ has compact support in $x$ in the interval $[a,b]$, then so does $h(\cdot,t)$ in nthe interval $[\eta(a,t),\eta(b,t)]$, for any $ t \in [0,T]$.
Moreover, defining 
\begin{equation} 
\label{1.51} 
E_{+}(t) = \int_{\eta(a,t)}^{\eta(b,t)} e^y h(y,t) \, dy \;\;\;\;\;\text{and}\;\;\;\;\
E_{-}(t)=  \int_{\eta(a,t)}^{\eta(b,t)} e^{-y} h(y,t) \, dy,
\end{equation}
 one has from \eqref{1.50} that
\begin{equation}
\label{1.52}
u(x,t)=\,e^{-|x|}\ast h(x,t)=\, e^{-x} E_{+}(t),\;\;\;\;\;\;\;x>\eta(b,t),
\end{equation}
and
\begin{equation}
\label{1.53}
u(x,t)=\,e^{-|x|}\ast h(x,t)=\, e^{x} E_{-}(t),\;\;\;\;\;\;\;x<\eta(a,t).
\end{equation}

Hence, it follows  that for $x>\eta(b,t)$ 
\begin{equation}
\label{1.54}
u(x,t) = -\p_x u(x,t)= \p_x^2 u(x,t)=e^{-x} E_{+}(t),
\end{equation}
and for $x<\eta(a,t)$ 
\begin{equation}
\label{1.55}
u(x,t) = \p_x u(x,t)= \p_x^2 u(x,t)= e^{x} E_{-}(t).
\end{equation}

Next, 
integration by part, \eqref{1.54}, \eqref{1.55}, and the equation in \eqref{1.1}  yield the identities
\begin{equation}
\label{1.56}
\begin{aligned}
E_{+}(0)&=\int_{-\infty}^{\infty} e^y h(y,0) dy = 
\int_{-\infty}^{\infty} e^y u_0(y) dy -\int_{-\infty}^{\infty} e^y \p_x^2 u_0(y) dy\\
& =\int_{-\infty}^{\infty} e^y u_0(y) dy +\int_{-\infty}^{\infty} e^y \p_x u_0(y) dy=0
\end{aligned}
\end{equation}
and  
\begin{equation}
\label{1.57}
\begin{aligned}
&\frac{dE_{+}(t)}{dt\;\;\;}\\
&=-\int_{-\infty}^{\infty} e^y u\p_x udy+
\int_{-\infty}^{\infty} e^y \p_x^2(u\p_x u)dy
-\int_{-\infty}^{\infty} e^y \p_x F(u) dy\\
&=e^y(\p_x(u\p_x u)-u\p_xu) |_{-\infty}^{\infty}
                                    -e^y\left(u^2+\dfrac{(\p_x u)^2}{2}\right)|_{-\infty}^{\infty} 
+\int_{-\infty}^{\infty} e^y \ F(u) dy\\
&=\int_{-\infty}^{\infty} e^y \left(u^2+\dfrac{(\p_x u)^2}{2}\right)(y,t) dy > 0.
\end{aligned}
\end{equation}
Therefore, in the life-span  of the solution $u(x,t)$,
$\,E_{+}(t)$ is an increasing function. Thus, from \eqref{1.56} it follows that $E_{+}(t)>0$ for $t\in(0,T]$.

Similarly, it is easy to see that $E_{-}(t)$ is decreasing with $E_{-}(0)=0$, therefore $E_{-}(t)<0$ for $t\in(0,T]$.

Taking $c_{\pm}(t)=E_{\pm}(t)$ we obtain \eqref{1.13c}.

\vskip.1in

Next, let us consider part (b). Since $ \,h(x,t) = (1 -\p_x^2) u(x,t)$ satisfies the equation 
\begin{equation}
\label{1.58}
\p_t h + u(x,t) \p_x h = -2 \p_xu(x,t) h,\;\;\;\;\;\;\;\;(x,t)\in\rr \times [0,T],
\end{equation}
an argument similar to that given in the proof of Theorem \ref{Theorem 2} shows that
\begin{equation}
\label{1.59}
\sup_{t\in [0,T]} \| h(t) e^{(1+\mu)|x|}\|_{\infty}\leq \tilde c \| h(0) e^{(1+\mu)|x|}\|_{\infty},
\end{equation}
with $\tilde c\,$ depending only on $M$ in \eqref{1.22} and $T$, and that for any $\,\theta\in (0,1)$
\begin{equation}
\label{1.60}
 \p_x^j u(t) \sim O(e^{-\theta |x|})\;\;\;\;\;\;\;\text{as}\;\;\;\;\;|x|\uparrow \infty\;\;\;\;\;\text{for}\;\;\;\;\;j=0,1,2.
\end{equation}

Thus, the definitions
in \eqref{1.51} make sense with the integrals extended  to  the whole real line
and the computations in \eqref{1.56}-\eqref{1.57}  can be carried out in the same fashion.
Finally, using \eqref{1.59} in \eqref{1.49} we obtain \eqref{1.13f}.

\end{proof} 

\vspace{3mm} \noindent{\large {\bf Acknowledgments}}
\vspace{3mm}\\ The authors would like to thank Prof. D. Holm for  useful comments concerning this work.
G. P.  was supported by an NSF grant. Y. Z. was supported by NSFC under grant no.
10501012 and Shanghai Rising-Star Program 05QMX1417.

\vskip.3in

\noindent{\bf A. Alexandrou Himonas}\\
Department of Mathematics\\
University of Notre Dame\\
Notre Dame, IN 46556\\
USA\\
E-mail: himonas.1@nd.edu\\

\noindent{\bf Gerard Misio\l ek}\\
Department of Mathematics\\
University of Notre Dame\\
Notre Dame, IN 46556\\
USA\\
E-mail: gmisiole@nd.edu\\

\noindent{\bf Gustavo Ponce}\\
Department of Mathematics\\
University of California\\
Santa Barbara, CA 93106\\
USA\\
E-mail: gmisiole@nd.edu\\

\noindent{\bf Yong Zhou}\\
Department of Mathematics\\
East China Normal University\\
Shangai 200062\\
China\\
E-mail: yzhou@math.ecnu.cn\\
and

\noindent Institute des Hautes \'Eudes Scientifiques\\
35, route de Chartres \\
F-91440 Bures-sur-Yvette\\
France\\


\begin{thebibliography}{EKPV1} 



\bibitem[BSS]{bss}
      Beals,  R.,  Sattinger, D., and  Szmigielski, J.,
{\em Multipeakons and the classical moment
problem},  Adv. Math.  {\bf 154}  (2000),  no. 2, pp. 229--257.


\bibitem[CH]{ch} 
Camassa, R. and Holm, D., 
\textit{An integrable shallow water equation with peaked solutions}, 
Phys. Rev. Lett. {\bf 71} (1993), pp. 1661-1664. 

\bibitem[Co]{C} 
Constantin, A.,
\textit{Finite propagation speed for the Camassa-Holm equation},
J. Math. Phys. {\bf 46} (2005), no 2, pp. 4

\bibitem[CoE]{CE} 
Constantin, A. and Escher, J.,
\textit{Global existence and blow-up for a shallow water equation},
Ann. Scuola Norm. Sup. Pisa Cl. Sci. {\bf 26}  (1998),  no. 2, pp. 303--328.

\bibitem[CoMc]{comc} 
Constantin, A. and McKean, H., 
\textit{A shallow water equation on the circle}, 
Comm. Pure Appl. Math. {\bf 52} (1999), pp. 949-982. 
 

\bibitem[CoS]{cs} 
Constantin, A. and Strauss, W., 
\textit{Stability of peakons}, 
Comm. Pure Appl. Math. {\bf 53} (2000), pp. 603-610. 

\bibitem[D]{d} 
Danchin, R., 
\textit{A few remarks on the Camassa-Holm equation}, 
Differential Integral Equations {\bf 14} (2001), pp. 953-988. 


\bibitem[DKT]{dkt} 
De Lellis, C., Kappeler, T. and Topalov, P., 
\textit{Low-regularity solutions of the periodic Camassa-Holm equation}, 
pre-print. 

\bibitem[EKPV1]{EKPV1} 
Escauriaza, L., Kenig, C. E., Ponce, G., and Vega, L., 
\textit{On unique continuation of solutions of Schr\"odinger  equations}, 
to appear in Comm. PDE 

\bibitem[EKPV2]{EKPV2} 
Escauriaza, L., Kenig, C. E., Ponce, G., and Vega, L., 
{\em On uniqueness properties of solutions of  the $k$-generalized KdV equations}, 
pre-print. 

\bibitem[FF]{ff} 
Fuchssteiner, B. and Fokas, A., 
\textit{Symplectic structures, their Backlund transformations and 
hereditary symmetries}, 
Phys. D {\bf 4} (1981/1982), pp. 47-66. 

\bibitem[F]{f} 
Fuchssteiner, B.,
\textit{Some tricks from the symmetry-toolbox for nonlinear equations:
generalization of the Camassa-Holm equation},
\textit{Physica D}  {\bf 95}, (1996), pp. 229-243.

\bibitem[He]{He} 
Henry, D., 
\textit{Compactly supported solutions of
the Camassa-Holm equation}, 
J. Nonlinear Math. Phys., {\bf  12} (2005), pp. 342-347. 

\bibitem[HM1]{hm1} 
Himonas, A. and Misio\l ek G., 
\textit{The Cauchy problem for an integrable shallow water equation},
Differential Integral Equations {\bf 14} (2001), pp. 821-831.


\bibitem[HM2]{hm2} 
Himonas, A. and Misio\l ek G., 
\textit{High-frequency smooth solutions and well-posedness 
of the Camassa-Holm equation}, 
Int. Math. Res. Not. {\bf  51} (2005), pp. 3135-3151. 

\bibitem[LO]{lo} 
Li, Y. and Olver, P., 
\textit{Well-posedness and blow-up solutions for an integrable nonlinearly 
dispersive model wave equation},  
J. Differential Equations {\bf 162} (2000), pp. 27-63. 

\bibitem[Mc]{mc} 
McKean, H., 
\textit{Breakdown of the Camassa-Holm equation}, 
Comm. Pure Appl. Math. {\bf 57} (2004), pp. 416-418. 

\bibitem[Mi]{mi} 
Misio\l ek, G., 
\textit{Classical solutions of the periodic Camassa-Holm equation},  
Geom. Funct. Anal. {\bf 12} (2002), pp. 1080-1104. 

\bibitem[Mo]{mo} 
Molinet, L., 
\textit{On well-posedness results for the Camassa-Holm equation 
on the line: A survey}, 
J. Nonlin. Math. Phys. {\bf 11} (2004), pp. 521-533.  

\bibitem[R]{rb} 
Rodriguez-Blanco, G., 
\textit{On the Cauchy problem for the Camassa-Holm equation}, 
Nonlinear Anal. {\bf 46} (2001), pp. 309-327. 

\bibitem[Z]{Z} Zhou, Y., 
{\em Infinite propagation speed for a shallow water equation}, 
pre-print.

\end{thebibliography}
\end{document}